\documentclass[oneside,12pt,a4paper,english]{amsart}
\usepackage[T1]{fontenc}
\usepackage[utf8]{inputenc}
\usepackage[english]{babel}
\usepackage{csquotes}
\usepackage{amsmath}
\usepackage{booktabs}
\usepackage{enumitem}
\usepackage{amssymb}
\usepackage{amsthm}
\usepackage{amstext}
\usepackage{pdflscape}
\usepackage{booktabs}
\usepackage{graphicx}
\usepackage{multirow}
\usepackage[table]{xcolor}
\usepackage{diagbox}
\usepackage{placeins}
\usepackage[most]{tcolorbox}
\usepackage{geometry}

\usepackage{hyperref}
\definecolor{Blueish}{HTML}{1c4197}

\newtcolorbox{Chaptersetting}[1][]{%
	colback=black!5,
	fonttitle=\bfseries,
	title={#1},
	coltitle=black,
	sharp corners,
	borderline west={2pt}{0pt}{Blueish},
	enhanced,
}

\theoremstyle{plain}

\newtheorem*{thm*}{Theorem}

\newtheorem*{lem*}{Lemma}

\theoremstyle{definition}

\usepackage{graphicx}

\begin{document}

	\title[The Mathematician's Assistant: Integrating AI into Research Practice]{The Mathematician's Assistant: Integrating AI into Research Practice}
	
	\author{Jonas Henkel}
	
	\address{\hspace{-7mm} 
Jonas Henkel,
Fachbereich Mathematik und Informatik,
Philipps-Universit\"at Marburg,
Campus Lahnberge,
35032 Marburg, Germany\newline
{\normalfont\ttfamily henkelj@mathematik.uni-marburg.de}}

\keywords{Artificial Intelligence, Mathematics, Large Language Models, Mathematical Reasoning, AI-Assisted Research, Benchmarking, MathArena, Open Proof Corpus, Gemini, GPT, Grok, AI Ethics, Human-AI Collaboration}

	\maketitle
    \begin{center}
   \today
    \end{center}
    
\begin{abstract}
	The rapid development of artificial intelligence (AI), marked by breakthroughs like \texttt{AlphaEvolve} and \texttt{Gemini Deep Think}, is beginning to offer powerful new tools that have the potential to significantly alter the research practice in many areas of mathematics. This paper explores the current landscape of publicly accessible large language models (LLMs) in a mathematical research context, based on developments up to August 2, 2025. Our analysis of recent benchmarks, such as MathArena and the Open Proof Corpus \cite{MathArena BalunovicEtAl25, Large Scale Study OPCDeEtAl25}, reveals a complex duality: while state-of-the-art models demonstrate strong abilities in solving problems and evaluating proofs, they also exhibit systematic flaws, including a lack of self-critique and a model depending discrepancy between final-answer accuracy and full-proof validity.

	Based on these findings, we propose a durable framework for integrating AI into the research workflow, centered on the principle of the augmented mathematician. In this model, the AI functions as a copilot under the critical guidance of the human researcher, an approach distilled into five guiding principles for effective and responsible use. We then systematically explore seven fundamental ways AI can be applied across the research lifecycle, from creativity and ideation to the final writing process, demonstrating how these principles translate into concrete practice. We conclude that the primary role of AI is currently augmentation rather than automation. This requires a new skill set focused on strategic prompting, critical verification, and methodological rigor in order to effectively use these powerful tools.
\end{abstract}
	\tableofcontents
\section{The New Ecosystem: AI as an Economic and Scientific Catalyst}
The development of artificial intelligence (AI) has seen exponential growth over recent years, profoundly influencing both economic landscapes and scientific advancement. In economic terms, AI has rapidly evolved into a substantial market sector. In 2024 alone, global private investment in AI reached an impressive 252 billion USD, with a notable 33.9 billion USD specifically allocated to generative AI technologies—an 8.5-fold increase compared to investments in 2022 \cite{StanfordHAI25}.

Beyond financial indicators, the practical adoption of AI within businesses has grown significantly. According to a comprehensive survey by McKinsey $\&$ Company, by 2024, 78~\% of globally surveyed businesses reported utilizing AI in their operations, up significantly from 55~\% in 2023, and merely 20~\% in 2017 \cite{StanfordHAI25}. This growing reliance highlights AI's escalating importance as a key driver of corporate innovation and operational efficiency.

Moreover, the impact of AI extends deeply into scientific research, marking historic milestones in fields traditionally dominated by human ingenuity.  Notably, the 2024 Nobel Prize in Chemistry recognized groundbreaking work in protein folding, awarding John Jumper and Demis Hassabis for their development of \texttt{AlphaFold} at Google DeepMind, alongside David Baker for his independent, pioneering research in computational protein design. \texttt{AlphaFold}, in particular, represents an AI-driven breakthrough that has fundamentally transformed biological research. Simultaneously, the 2024 Nobel Prize in Physics was awarded for pioneering research in AI fundamentals, specifically recognizing John Hopfield and Geoffrey Hinton for their seminal contributions to artificial neural networks and machine learning.

This paper does not address broader societal risks associated with AI, including algorithmic bias, environmental impact, and economic displacement. The author acknowledges the importance of these issues. The scope of this article is to provide a practical framework for the effective and responsible use of AI within the mathematical research workflow, so that everybody can make an informed decision whether and how to use AI tools. %, it's important to understand that AI is also a powerful tool. To make an informed decision about its application in research, we need to understand its practical use. Therefore, 

\section{Landmark Achievements in AI-Driven Mathematical Discovery}\label{sec: Landmark Achievements in AI-Driven Mathematical Discovery}

The influence of artificial intelligence on mathematics is becoming increasingly evident, though its impact must be carefully put in context. Although significant breakthroughs in core theoretical fields, such as resolving major conjectures or winning a Fields Medal, remain a distant prospect, AI's progress on difficult mathematical problems has accelerated significantly. This development is best understood not as an automation of mathematical discovery, but as the emergence of a powerful new class of problem-solving tools. Their current strengths are most clearly demonstrated in two distinct areas: high-level competitions that test for creative reasoning, and complex optimization problems that are foundational to computational science.

In 2024, AI-based systems \texttt{AlphaGeometry} and \texttt{AlphaProof} required human experts to translate International Mathematical Olympiad (IMO) problems into machine-readable code, and computations spanned multiple days, ultimately securing a silver medal. Only one year later, in 2025, Google's advanced model, \texttt{Gemini Deep Think}, marked a substantial leap forward. Officially permitted to compete autonomously in the IMO, \texttt{Gemini Deep Think} secured a gold medal for the first time, solving problems entirely autonomously in natural language within the stringent 4.5-hour competition time limit \cite{ReutersCaiSingh25, Google IMO DeepMindLuongLockhart25}.

Let us have a closer look at this achievement. The formal mathematical knowledge required for the IMO is pre-university, covering topics such as geometry, number theory, and algebra. The problems are notorious for their extreme difficulty, demanding not just knowledge, but a profound level of ingenuity, creative synthesis, and logical rigor that far exceeds standard coursework. Be it students or professional mathematicians, taking part in these challenges requires significant preparation to compete. Therefore, the IMO serves as an exceptional benchmark not for accumulated knowledge, but for genuine, creative problem-solving ability—a skill long considered a key feature of human intellect. The success of \texttt{Gemini Deep Think} is a landmark achievement because it demonstrates that an AI can operate at this high level of creative reasoning, autonomously and under the same time restrictions as its human competitors. The full set of solutions provided by the model is publicly available for analysis \cite{DeepMindIMO25}.

\begin{figure}[b]
	\centering
	\includegraphics[width=0.7\linewidth]{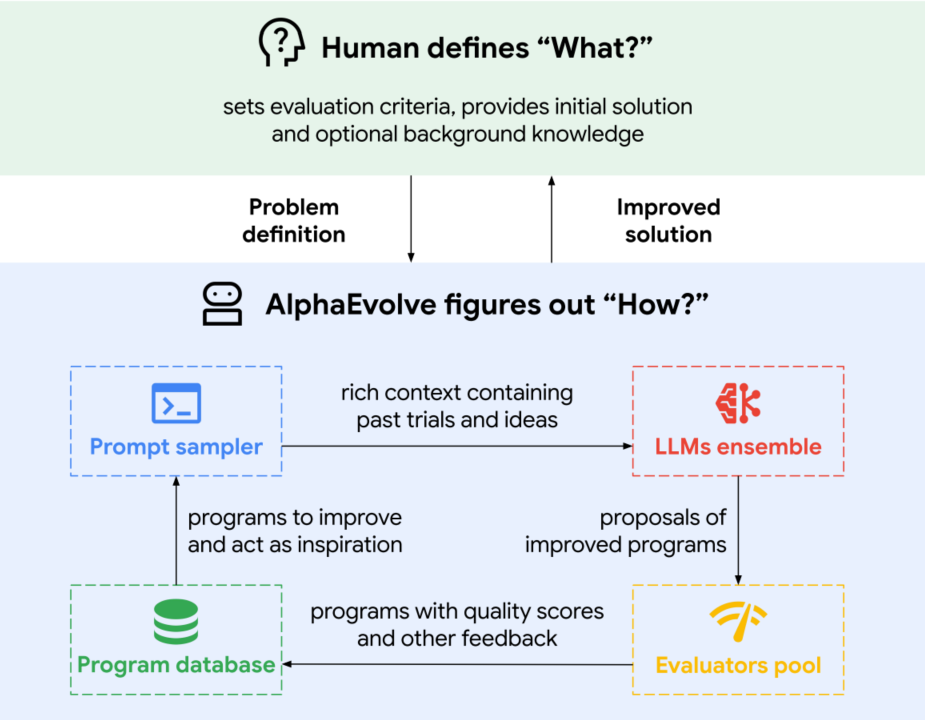}
	\caption{Operating Principle of \texttt{AlphaEvolve} \cite[Fig. 1]{AlphaEvolve BalogEtAl25}}
	\label{fig:alpha-evolve}
\end{figure}
Parallel to these advances in creative reasoning, AI has demonstrated profound capabilities in a different domain: complex optimization. Google DeepMind further pushed these boundaries with \texttt{AlphaEvolve} \cite{Google AlphaEvolve DeepMind25}, an advanced AI model designed specifically for internal research. After an initial setup by a human researcher, \texttt{AlphaEvolve} operates autonomously, leveraging various Gemini language models in an evolutionary process to generate and rigorously test new code variants, iteratively refining successful solutions to continually enhance performance (see Figure \ref{fig:alpha-evolve}).
In a preprint, the developers describe applying \texttt{AlphaEvolve} to a set of over 50 challenging problems from fields like analysis, number theory, geometry, and combinatorics, many of which involve finding optimal algorithms or constructions \cite{AlphaEvolve BalogEtAl25}. The problems are all clearly defined, and the background knowledge required to understand them is that of an undergraduate student. The results were notable: in 20~\% of the cases, the model discovered novel solutions that surpassed the best-known existing results, thereby pushing the state of the art in these specialized domains. A specific example from the field of analysis is the \enquote{second autocorrelation inequality}. For a non-negative function $f \colon \mathbb{R} \to \mathbb{R}$, the task is to find the smallest constant $C_2$ for which the inequality
\[
\|f * f\|_{2}^{2} \le C_2 \|f * f\|_{1} \|f * f\|_{\infty}
\]
holds. While it was known that $0.88922 \le C_2 \le 1$, \texttt{AlphaEvolve} discovered a new function that established a slightly better lower bound of $0.8962 \le C_2$, see \cite[B.2.]{AlphaEvolve BalogEtAl25}.

One notable breakthrough occurred in computational mathematics: \texttt{AlphaEvolve} discovered an innovative algorithm for multiplying two $4\times4$ matrices with only 48 scalar multiplications. This advancement marks the first improvement in over 56 years, surpassing the previously best-known Strassen algorithm, which required 49 multiplications. Such discoveries are of immense practical importance, as more efficient algorithms for core operations like matrix multiplication are crucial for accelerating the training and inference of large-scale AI systems themselves. Further specific examples include significant advancements in Paul Erdős's \enquote{Minimum Overlap Problem} and the \enquote{Kissing Number Problem} in 11 dimensions.

Taken together, the IMO results and the achievements of \texttt{AlphaEvolve} illustrate the two primary fronts of AI's advance in mathematics: creative, proof-based reasoning on one hand, and autonomous, large-scale algorithmic optimization on the other.
\section{From Frontier Models to Publicly Affordable Tools}\label{sec:From Frontier Models to Publicly Available Tools}
The landmark achievements of specialized systems like \texttt{Gemini Deep Think} and \linebreak\texttt{AlphaEvolve} are powerful demonstrations of the frontier of AI's mathematical capabilities. However, these frontier models rely on computational resources and are available through access schemes—such as internal-only use (\texttt{AlphaEvolve}) or premium-tier subscriptions (\texttt{Gemini Deep Think}) —that place them far beyond the standard toolkit of most researchers. Their performance is therefore an indicator of future potential rather than a reflection of the more widely affordable tools that are the focus of this paper. 

To bridge this gap and assess the practical utility of AI for the everyday researcher, it is essential to evaluate the performance of widely accessible large language models (LLMs). This requires standardized, rigorous, and uncontaminated benchmarks that can reliably measure the mathematical reasoning abilities of these widely available tools. The analysis is grounded in the state of the AI landscape as of August 2, 2025, placing it immediately before the release of OpenAI's \texttt{GPT 5} on August 7. This rapid pace of development is a central theme of this paper. Since the \texttt{GPT 5} update consolidated the benchmarked \texttt{o}-series into a new system for which no independent evaluations yet exist, our findings establish the critical mid-2025 baseline against which this new generation of models will be measured.

\subsection{MathArena: Answer-Based and Proof-Based Performance}
A fundamental challenge in evaluating the mathematical reasoning capabilities of LLMs lies in ensuring the novelty of test problems. It is often difficult to determine whether a model's successful performance is due to its genuine problem-solving abilities or its ability to recall solutions encountered during its extensive training phase. This phenomenon is known as data contamination.

To address this issue, benchmarks such as MathArena have been developed to provide a standardized and rigorous assessment framework. The methodology of MathArena is based on a key insight: recurring, high-stakes mathematical competitions serve as a continuous source of novel, high-quality problems. By exclusively sourcing tasks from recent events like the American Invitational Mathematics Examination (AIME) and the Harvard-MIT Mathematics Tournament (HMMT), this approach ensures that models are evaluated on genuinely unseen material.

\subsubsection{Answer-Based Performance}
Between January and May 2025, MathArena evaluated 30 LLMs on problems sourced from four prominent competitions: the American Invitational Mathematics Examination (AIME), the Harvard-MIT Mathematics Tournament (HMMT), the Brown University Math Olympiad (BRUMO), and the Stanford Math Tournament (SMT) \cite{MathArena BalunovicEtAl25}. It is crucial to contextualize these contests. While all require creative reasoning, their problems are significantly less demanding than those of the International Mathematical Olympiad (IMO). The AIME, for instance, serves as a qualifying stage for the USAMO, which in turn is part of the selection process for the US IMO team. The human participants are typically exceptionally talented pre-university students, and for this portion of the benchmark, an answer-based evaluation was employed, where correctness was determined solely by the final numerical answer.

\begin{figure}[h!]
	\centering
	\includegraphics[width=1\linewidth]{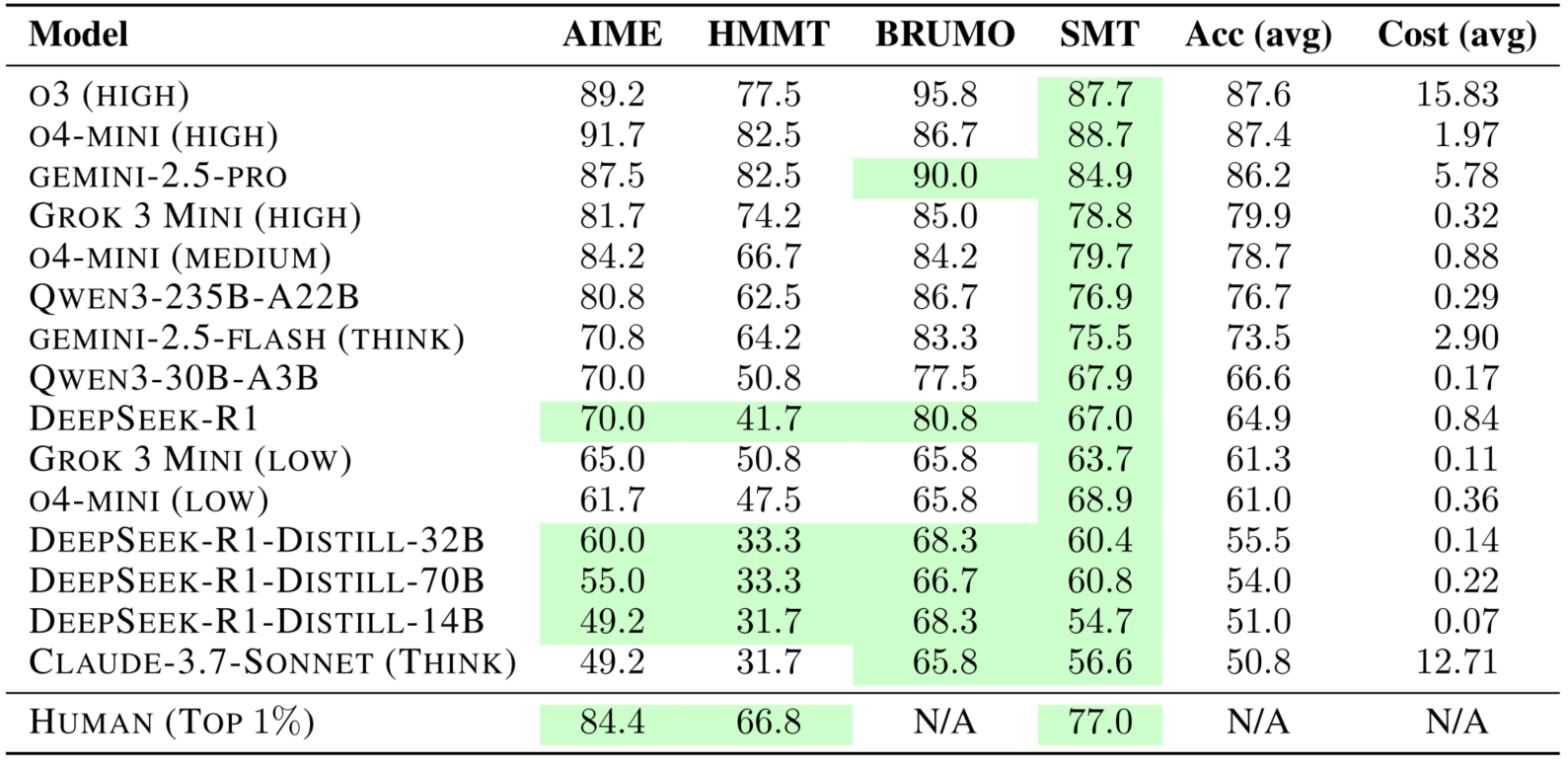}
	\caption{MathArena: Main results on numerical answer evaluation on competitions between January-May 2025. Green cells denote models released after the respective competition date \cite[Tab. 2]{MathArena BalunovicEtAl25}.}
	\label{fig:benchmarks}
\end{figure}

As shown in Figure \ref{fig:benchmarks}, models from Google (\texttt{Gemini 2.5 Pro Exp}) and OpenAI (\texttt{o3} and \texttt{o4 mini high}) dominated the tests, achieving scores that surpassed the performance of the top 1$~\%$ of human participants. All problems and model-generated solutions are publicly available on the benchmark's website \cite{MathArena Website25}.

\subsubsection{Proof-Based Performance}
The assessment of LLM capabilities becomes more nuanced when shifting to proof-based evaluations. This format requires a definitive answer and a complete, step-by-step logical reasoning process. This form of human-led evaluation is central to benchmarks sourced from elite competitions such as the United States of America Mathematical Olympiad (USAMO) and the International Mathematical Olympiad (IMO). Crucially, these competitions are also defined by a significant escalation in problem difficulty. The USAMO, for instance, is an invitational event where qualification is largely determined by top performance in preceding contests like the AIME. Some of these problems are designed to challenge the world's most elite pre-university talent, a standard that is formidable even for professional mathematicians.

Against this backdrop, the performance of publicly affordable LLMs provides a clear baseline. At the 2025 IMO, the top-performing accessible model, \texttt{Gemini 2.5 Pro}, scored  31.55~\%. In their analysis, the authors highlighted the considerable gap that remains to the threshold for a human medal-winning performance \cite{MathArena IMO25 authors}. A similar capability was demonstrated on the USAMO, where leading models also achieved scores of around 30~\% \cite{USAMO PetrovEtAl25}.

These results demonstrate the boundary of their reasoning skills. The scores do not indicate poor logic; we should also consider them in light of how difficult the problems were. This performance delineates the frontier of publicly affordable AI and reveals a significant gap when compared to the near-perfect scores of specialized, state-of-the-art models like \texttt{Gemini Deep Think}.

The performance decrease from pre-olympiad to the olympiad-style format is due to an increase in problem difficulty and a change in the evaluation method. This raises a critical question: How much of the performance drop is attributable to the intrinsic difficulty of the problems, versus the specific challenge of generating a coherent and logically sound proof? To isolate this variable and diagnose the failure modes in proof generation, a systematic, large-scale analysis of the proofs themselves became necessary.
\subsection{The Open Proof Corpus: A Deeper Dive into Proof Validity}
This need was addressed by the introduction of the Open Proof Corpus (OPC) on June 23, 2025. The OPC is a large-scale, human-evaluated dataset comprising over 5,000 LLM-generated proofs. The problems are sourced from a diverse mix of prestigious competitions, ranging from national contests like the British Mathematical Olympiad to international events such as the IMO and undergraduate-level challenges like the Putnam Competition. Its main goal is to analyze the proof-generation capabilities of LLMs in more detail, moving beyond final-answer accuracy to check the validity of the reasoning itself \cite[p. 1, 5]{Large Scale Study OPCDeEtAl25}.

\subsubsection{Overall Model Performance in Proof Generation}
\begin{figure}[h!]
		\centering
		\includegraphics[width=0.7\linewidth]{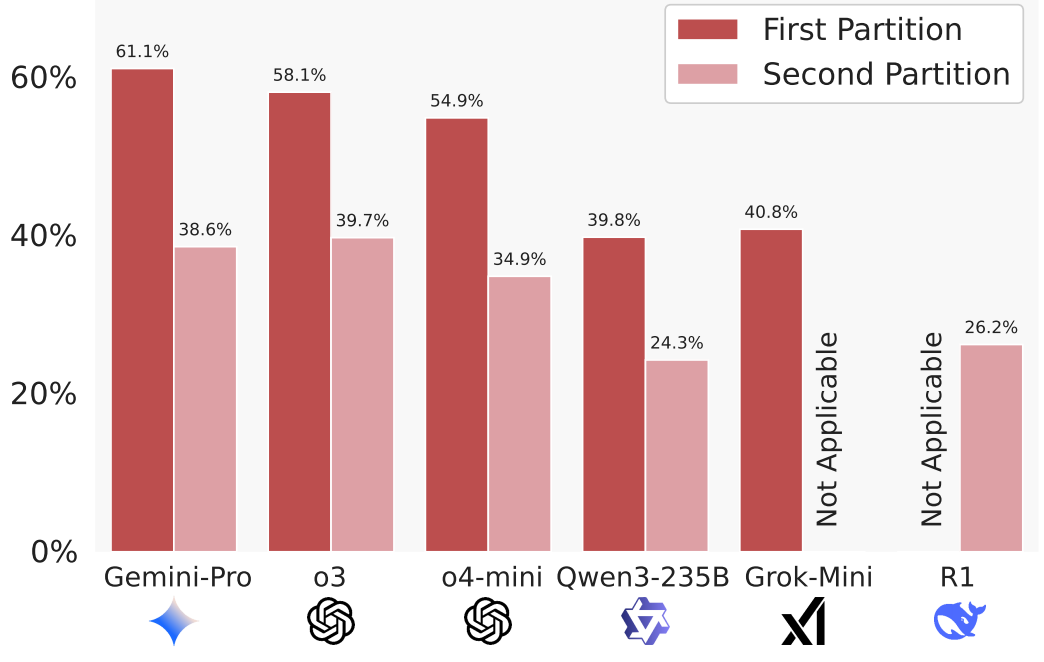}
		\caption{Average proof correctness of various models on the OPC. Data is split into two partitions. First, resp. second was answered by all models except for \texttt{DeepSeek-R1}, resp. \texttt{Grok 3 mini} \cite[Fig. 3]{Large Scale Study OPCDeEtAl25}.}
		\label{fig:proof-benchmarks}
	\end{figure}

The OPC analysis reveals that leading LLMs demonstrate a commendable, although lower, performance compared to their results on purely numerical evaluations. To ensure fair comparisons among models released at different times, the study partitioned the dataset. Specifically, the model \texttt{DeepSeek-R1} was introduced mid-way through the project, replacing \texttt{Grok 3 mini}. The results, shown in Figure \ref{fig:proof-benchmarks}, consistently identify \texttt{Gemini 2.5 Pro} and \texttt{o3} as the top performers across both partitions. The general performance drop observed in the second partition is attributed to the inclusion of a more challenging set of problems. Across the entire corpus, a notable 43~\% of all generated proofs were deemed correct by human evaluators \cite[p. 6]{Large Scale Study OPCDeEtAl25}.

	\subsubsection{Discrepancy Between Final-Answer Accuracy and Proof Validity}
\begin{figure}[t]
		\centering
		\includegraphics[width=0.7\linewidth]{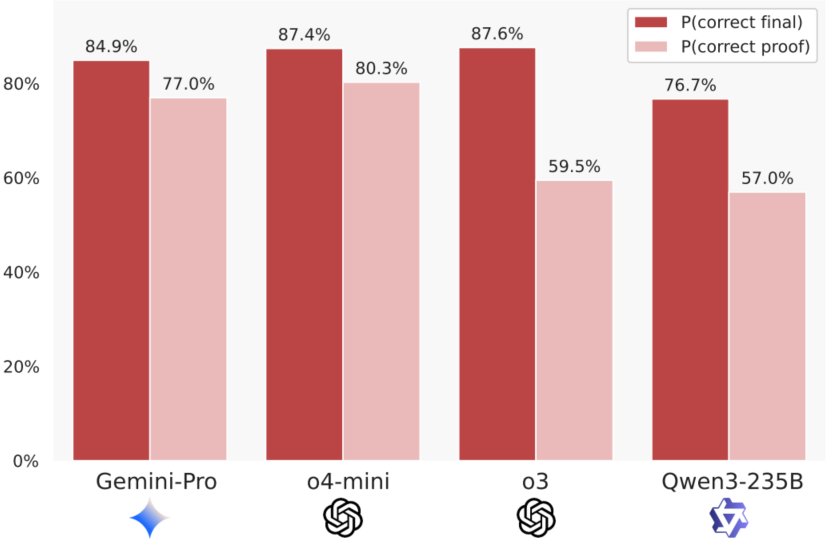}
		\caption{Comparison of final-answer accuracy and proof correctness on the MathArena subset \cite[Fig. 5]{Large Scale Study OPCDeEtAl25}}
		\label{fig:proof-vs-eval-benchmark}
	\end{figure}

	Using the MathArena subset, the OPC study empirically confirms that final-answer accuracy is a poor proxy for proof quality. The extent of this discrepancy, however, is highly model-dependent. As illustrated in Figure \ref{fig:proof-vs-eval-benchmark}, models like \texttt{Gemini 2.5 Pro} exhibit only a modest performance drop of approximately 8~\% when full proof validity is required. In stark contrast, other high-performing models, such as \texttt{o3}, experience a dramatic decline of nearly 30~\%. This finding is particularly insightful, as it suggests that LLMs can possess a strong, correct \enquote{intuition} for a final result, even when their underlying logical justification is flawed.

\subsubsection{LLMs as Evaluators of Mathematical Proofs}
\begin{figure}[b]
		\centering
		\includegraphics[width=0.7\linewidth]{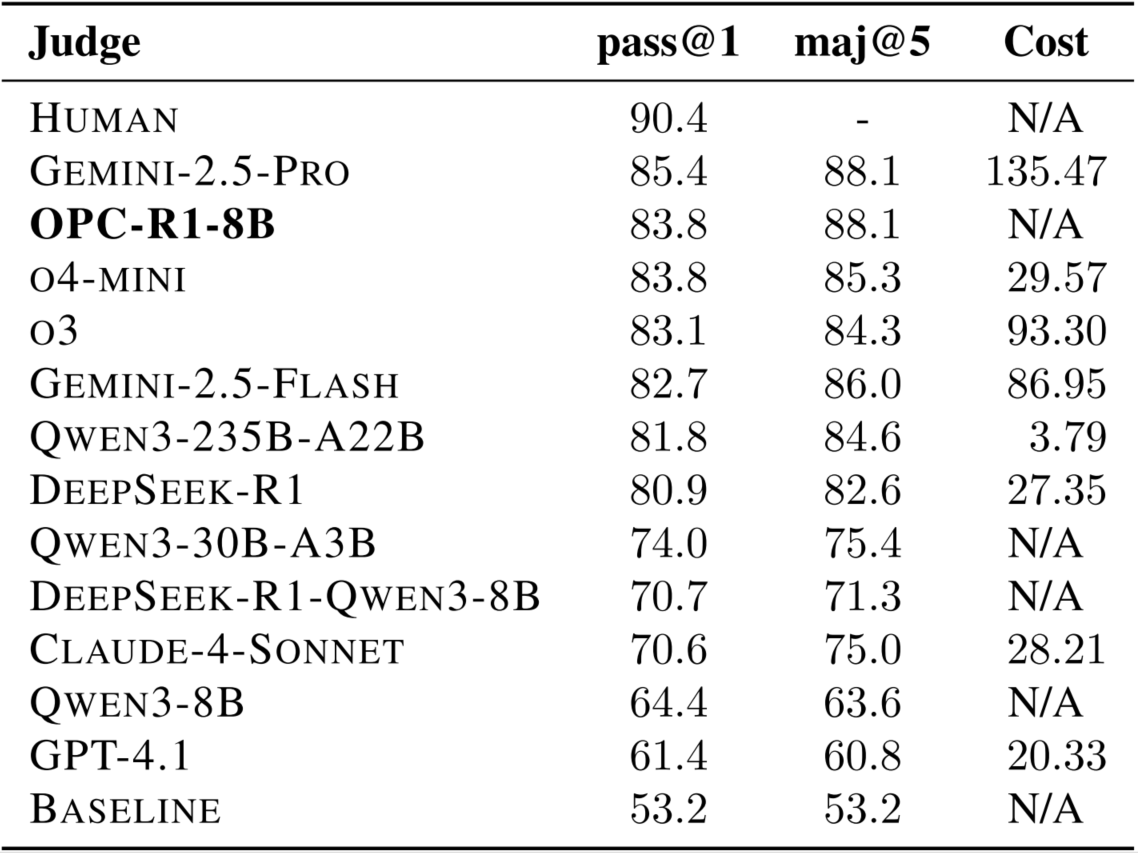}
		\caption{Benchmarking LLMs as proof graders. Costs are given in USD \cite[Tab. 2]{Large Scale Study OPCDeEtAl25}}
		\label{fig:j1}
	\end{figure}

	A surprising result from the OPC study is the proficiency of LLMs as evaluators, or \emph{judges,} of mathematical proofs. As shown in Figure \ref{fig:j1}, the leading models achieve a judgment accuracy that approaches human-level performance (e.g., \texttt{Gemini 2.5 Pro} at 85.4~\% compared to a human baseline of 90.4~\%). 
	
\begin{figure}[b]
		\centering
		\includegraphics[width=0.7\linewidth]{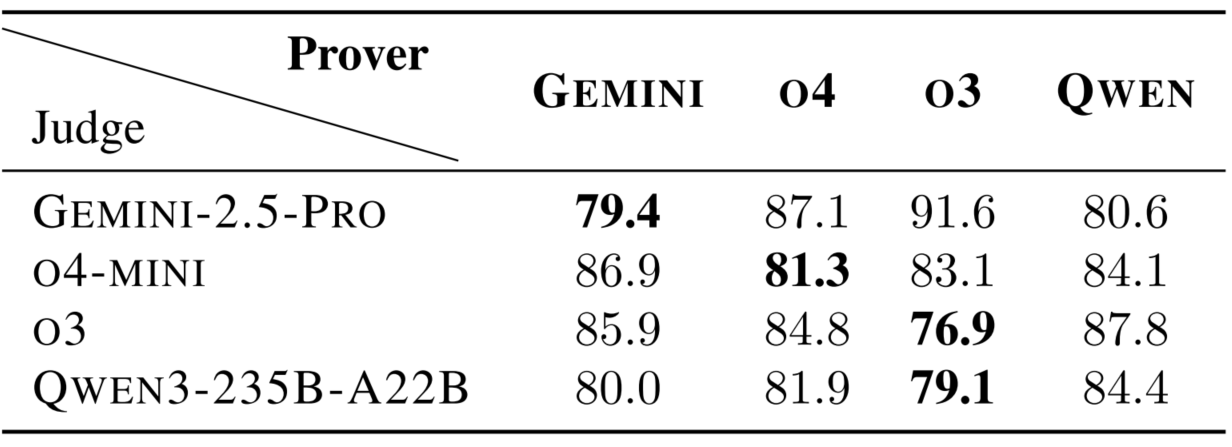}
		\caption{Judgement accuracy breakdown split by solver, highlighting the lowest score for each judge \cite[Tab. 3]{Large Scale Study OPCDeEtAl25}.}
		\label{fig:j2}
	\end{figure}

	However, a deeper analysis, presented in Figure \ref{fig:j2}, reveals a critical weakness. The table breaks down the judgment accuracy of each LLM when evaluating proofs generated by various other models. The boldface entries, which represent the lowest score for each judge, consistently lie along the main diagonal. This indicates that models are significantly worse at identifying errors in proofs generated by themselves. This \emph{self-critique blindness} is a crucial limitation, suggesting that for iterative refinement tasks, employing a combination of different models for generation and verification may be a more robust strategy.

	
	\subsubsection{Common Failure Modes in Proof Generation}\label{subsubsec: Common Failure Modes in Proof Generation}
	The analysis of incorrect proofs within the OPC reveals several recurring failure modes common to many LLMs. These include:
	\begin{itemize}
		\item \textbf{Overgeneralization:} A frequent error where a model validates a statement for a specific or simplified case and then incorrectly concludes its general validity without providing a rigorous, comprehensive proof.
		\item \textbf{Flawed Logical Steps:} Models often exhibit flawed reasoning, particularly when manipulating inequalities or applying geometric properties. This can lead to logically unsound arguments.
		\item \textbf{Reluctance to Admit Failure:} A notable behavioral trait is the models' tendency to produce a faulty proof rather than explicitly stating an inability to solve the problem. This occurs even when they are prompted to indicate uncertainty, highlighting a fundamental limitation in their self-awareness capabilities \cite[p. 6, Appendix C, p. 16-17]{Large Scale Study OPCDeEtAl25}.
	\end{itemize}
	
	\subsubsection{The Impact of Best-of-$n$ Sampling on Proof Quality}\label{subsubsec: The Impact of Best-of-$n$ Sampling on Proof Quality}
	The study also demonstrates that proof generation performance can be significantly enhanced through \enquote{best-of-$n$} sampling, a strategy that involves generating multiple candidate solutions and selecting the most promising one. For the model \texttt{o4 mini}, for instance, employing ranking-based selection methods dramatically improves the success rate from 26~\% (for a single attempt, or pass@1) to 47~\% when selecting the best out of eight generated proofs. This underscores the substantial potential of advanced sampling and selection techniques to boost the effective reasoning capabilities of LLMs \cite[p. 8, Fig. 6 \& 7]{Large Scale Study OPCDeEtAl25}.
\subsection{Competitions in July 2025}\label{subsec:Competitions in July 2025}
The competitive landscape continued to evolve rapidly into mid-2025. The release of \texttt{Grok 4} by xAI on July 10, 2025, introduced a new major contender, which was quickly reflected in subsequent benchmarks. Prior to its release, competitions were largely dominated by Google and OpenAI.

\texttt{Grok 4} quickly established itself at the top of the MathArena leaderboards. Averaging 89~\% correctness across several competitions within the benchmark, it secured the first-place ranking, narrowly outperforming the second-place model \texttt{o4 mini high} (87~\%) \cite{MathArena Website25}. Together with the recent release of \texttt{GPT 5}, these results point to a highly dynamic and competitive field among accessible LLMs, with new architectures rapidly challenging the established market leaders.
\subsection{FrontierMath: A Benchmark at the Research Frontier}
\begin{figure}[b]
	\centering
	\includegraphics[width=\linewidth]{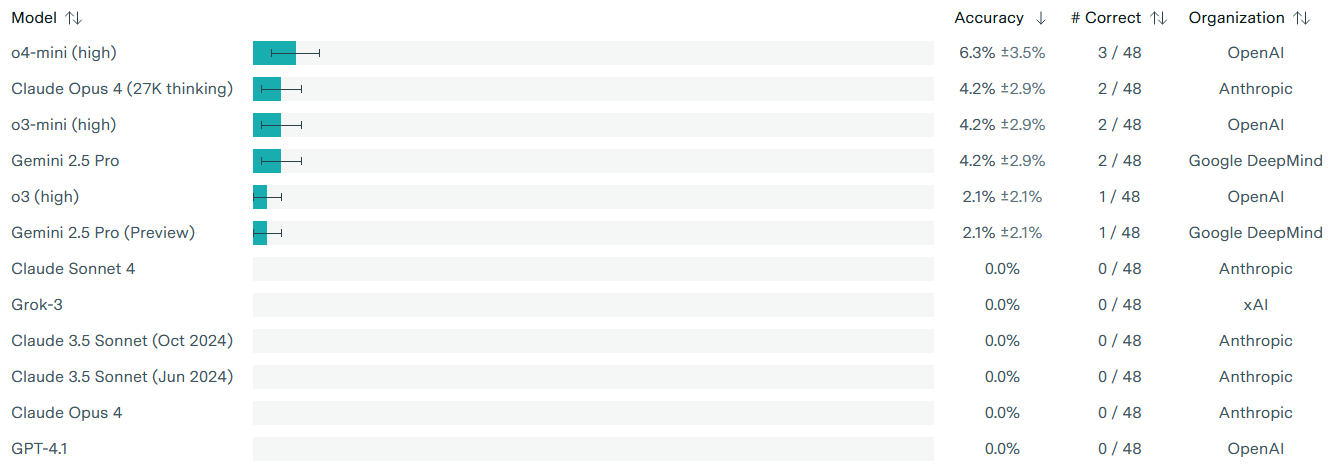}
	\caption{Performance of leading LLMs on the 50 problems of FrontierMath's Tier 4, the most difficult subset of the benchmark \cite{FrontierMath overview}.}
	\label{fig:frontiermath-leaderboard}
\end{figure}
To evaluate AI abilities on solving problems characteristic of modern mathematical research, the FrontierMath benchmark was introduced in late 2024 \cite{FrontierMath BarkleyEtAl24}. It comprises 350 original problems, structured into four tiers of increasing difficulty. The most challenging of these is Tier 4, an expansion of 50 \enquote{extremely difficult} problems developed by mathematics professors and postdoctoral researchers. However, the benchmark's creation was overshadowed by a significant controversy. As Epoch AI, the benchmark's creator, later clarified, the project was commissioned by OpenAI, which in return received ownership of most problems and their solutions. This arrangement was not transparently communicated to the contributing mathematicians from the outset \cite{EpochAIClarification25, HarrisFrontierMathScandal25}.

The FrontierMath controversy highlights the importance of examining AI benchmarks' funding, data access agreements, and overall independence. It illustrates the imbalance of power between big tech companies and individuals.  Due to the lack of transparency and potential data contamination, it is difficult to independently verify the performance claims of OpenAI, the funding partner. In order to maintain a conservative and scientifically rigorous analysis, we will set aside the results of the OpenAI models for this particular benchmark. Instead, we focus on the performance of the next-leading models from other developers on the publicly reported Tier 4 results. Among these, Google's \texttt{Gemini 2.5 Pro} and Anthropic's \texttt{Claude Opus 4 (27K thinking)} both achieved a success rate of 4.2~\% on the Tier 4 problems (see Figure \ref{fig:frontiermath-leaderboard}).

While this is a low absolute score, it represents a significant improvement over the performance of under 2~\% reported for the older models in the original arXiv publication \cite{FrontierMath BarkleyEtAl24}. These performances must be seen relative to the immense difficulty of the tasks: Ken Ono, a professor at the University of Virginia who collaborated on stating problems, remarked: \enquote{We designed problems that are our best guesses for challenges that would overwhelm AI. I can barely solve some of the problems from my area of expertise, so I hope AI entities don't get any of them right} \cite{FrontierMath About}. So this low score directly reflects the immense difficulty of the problems - an assessment shared by Fields medalist Terence Tao, who explained:
\begin{quote}
	\enquote{These are extremely challenging. I think that in the near term basically the only way to solve them, short of having a real domain expert in the area, is by a combination of a semi-expert like a graduate student in a related field, maybe paired with some combination of a modern AI and lots of other algebra packages…} \cite{FrontierMath About}
\end{quote}
This expert view underscores that the most effective path forward lies in a collaborative framework, leveraging the distinct strengths of both human intellect and artificial intelligence.

\section{The AI-Augmented Mathematician: A Practical Guide}
The preceding analysis demonstrates both the immense potential and the inherent limitations of modern AI models. This raises the critical question of how a researcher can responsibly navigate this new landscape. In 2023, Fields Medalist Terence Tao envisioned a future, predicting that by 2026, AI could serve as a \enquote{trustworthy co-author}. However, he included a crucial condition: this would only be possible when the technology is \enquote{used properly} \cite{TaoAIAnthology23}. This section is dedicated to defining what \enquote{proper use} entails for a working mathematician. To harness the power of AI effectively while mitigating its risks, the modern mathematician must adopt a new, strategic approach to their work. This guide outlines the core principles, essential tools, and practical applications for integrating AI into the mathematical workflow.
\subsection{Guiding Principles for AI Collaboration}\label{sub: Guiding Principles for AI Collaboration}
The integration of AI into the daily workflow of a mathematician is not merely about using a new software, but about adopting a new collaborative mindset. The models are powerful, but their effective and safe use depends on adhering to several core principles:

\begin{enumerate}
	\item \textbf{The Copilot, Not the Pilot:} The most fundamental principle is to understand that AI is a tool, a \emph{co-pilot} that assists, suggests, and computes. The mathematician remains the \emph{pilot}, responsible for direction, critical judgment, and the final verification of all results. This role includes the crucial task of determining when to use AI to accelerate progress and when doing so might be a time-consuming distraction, thus maintaining a productive balance between human intellect and machine assistance.
	
	\item \textbf{The Principle of Critical Verification:} No output from an LLM, whether it is a proof, a calculation, or a literature summary, should be accepted without rigorous verification. Every claim must be critically evaluated and, where possible, cross-checked against established literature or independent computation.
	
	\item \textbf{Understanding the Non-Human Nature of AI:} It is dangerous to anthropomorphize AI. Models do not \enquote{understand} in a human sense and do not forget incorrect conclusions, even after being corrected. This persistent error propagation requires constant vigilance and, at times, starting new sessions to ensure a clean reasoning slate.
	
	\item \textbf{The Art of Prompting and Model Selection:} Effective collaboration requires skill. This includes the art of asking the right questions (\emph{prompting}) and the experience to choose the right model for the right task (e.g., reasoning vs. writing vs. translation). The model landscape evolves rapidly, demanding continuous adaptation from the user.
	
	\item \textbf{The Experimental Mindset:} Leveraging AI to its full potential requires curiosity and a willingness to experiment. This includes testing different models, exploring creative prompts, and pushing the boundaries of what the tools can achieve in one's specific research area.
\end{enumerate}
With these principles in mind, let's turn to the concrete tools available to the modern mathematician.

\subsection{A Survey of Current AI Tools for the Working Mathematician}

The current AI landscape offers a diverse range of tools, from versatile, large-scale AI assistants to highly specialized applications. To provide a clear comparison, we first focus on the three major conversational AI platforms that are most relevant for complex mathematical tasks.

\subsubsection{Google DeepMind's Portfolio}\label{subsubsec: Google DeepMind's Portfolio}
\begin{table}[h!]
	\centering
	\caption{Overview of Google DeepMind Models}
	\label{tab:google_models}
	\begin{tabular}{p{3.5cm} p{3.5cm} p{6cm}}
		\toprule
		\textbf{Model} & \textbf{Availability} & \textbf{Key Strengths} \\
		\midrule
		\texttt{Gemini Deep Think} & available via AI Ultra €275/month & Frontier scientific reasoning; won IMO gold medal. \\
		\midrule
		\texttt{Gemini 2.5 Pro} & Free via AI Studio & Dominates Benchmarks; powerful multimodal capabilities (text, image, audio), large context window, creativity bar. \\
		\midrule
		\texttt{Deep Research} & Free (15 queries/month) & Specialized for academic research; optimized access to scholarly data. \\
		\bottomrule
	\end{tabular}
\end{table}

Google DeepMind's offerings reflect a dual strategy: pushing the scientific frontier with state-of-the-art models while simultaneously integrating powerful AI into a widely accessible ecosystem. This approach manifests in two distinct tiers: elite, high-cost access to frontier models for professionals, and free, powerful tools for the general public.

At the top of this strategy is the Google AI Ultra subscription. Priced at €274.99 per month, this premium tier provides access to Google's most advanced capabilities, most notably \texttt{Gemini Deep Think}. This is the model that achieved a gold medal at the International Mathematical Olympiad, representing the pinnacle of Google's mathematical reasoning abilities. Its availability through a high-cost subscription tier marks a clear distinction between elite frontier models and the powerful tools available through standard, more affordable subscription plans.

In contrast, Google makes its highly capable \texttt{Gemini 2.5 Pro} model, which has consistently dominated public benchmarks (see Section \ref{sec:From Frontier Models to Publicly Available Tools}) available for free through the \emph{Google AI Studio} platform (\url{https://aistudio.google.com}). This tool is particularly relevant for mathematicians due to several key features:
\begin{itemize}
	\item \textbf{Massive Context Window:} Each chat session supports a context window of 1 million tokens, equivalent to roughly 750,000 words. This is substantially larger than the context windows offered by competitors like OpenAI and xAI, and allows users to upload and analyze multiple large documents, such as textbooks or lengthy research papers, within a single, coherent conversation.
	
	\item \textbf{Demonstrated Proof-Generation Prowess:} Its strong performance in mathematical reasoning is empirically validated by the Open Proof Corpus (OPC) study. Specifically, it achieved the highest score for proof correctness (see Figure \ref{fig:proof-benchmarks}). It also demonstrated a remarkable consistency between its intuition and logical rigor, with its proof correctness on the MathArena subset dropping by only about 8 percentage points from its final-answer accuracy (Figure \ref{fig:proof-vs-eval-benchmark}). Furthermore, the model excelled as an evaluator, achieving the highest average accuracy when judging proofs (Figure \ref{fig:j1}).
	
	\item \textbf{Adjustable Creativity and Output Control:} A crucial feature for mathematical exploration is the adjustable \emph{temperature parameter}, which controls the model's creativity. A low temperature enforces deterministic, factual responses suitable for literature analysis, while a high temperature encourages novel solution paths, making it a powerful tool for generating new research ideas, albeit at an increased risk of hallucination. Users can also adjust the output length and the model's \enquote{thinking} time to further tailor the response.

	\item \textbf{Persistent Memory:} It is vital to understand that the model's memory is persistent within a session. Besides the advantages this causes, it also leads to unwanted behavior: it will not \enquote{forget} incorrect conclusions even if corrected and it may use those conclusions again. To ensure a clean reasoning, starting a new chat is sometimes necessary. A recommended best practice is to dedicate a single chat session to a specific problem, effectively creating a specialized assistant for that task.
	
\end{itemize}
The model is capable of writing and executing code directly within the interface. However, for maintaining transparency and control, it is recommended to execute the provided code manually, see Subsection \ref{subsubsec: A note on coding}. The platform also offers an integrated web search, though it is important to note that enabling this feature can sometimes degrade the quality of pure mathematical reasoning.\\
Alongside the versatile AI Studio, Google provides a specialized tool for in-depth literature synthesis: \texttt{Gemini Deep Research}. This tool is accessible for free with a limit of 15 queries per month (\url{https://gemini.google/overview/deep-research/}). A single query initiates an automated process. Before executing the full 5-10 minute search, the tool first generates a research plan, which the user can then review and modify—allowing them to correct, expand, or shorten the scope of the inquiry. Based on this finalized plan, the AI performs iterative web searches, analyzes documents, and uses its findings to refine subsequent searches. The final output is a synthesized summary with verifiable source citations, making it ideal for quickly getting up to speed on a new research topic. The collected statements can be verified by opening the respective links or documents. Table \nolinebreak \ref{tab:google_models} provides a concise summary of Google DeepMind's offerings.

\subsubsection{OpenAI's Product Stack}

OpenAI's ecosystem is characterized by a diverse suite of models and tools, each tailored to specific tasks. This specialization makes it crucial for the user to understand the distinct capabilities of each component to leverage the platform effectively. The following analysis describes the product stack as it existed for most of 2025, which formed the basis for the benchmark results discussed in this paper. An overview of the primary models and their functions within the different subscription tiers of that era is provided in Table \ref{tab:openai_models}.

As established in Section \ref{sec:From Frontier Models to Publicly Available Tools}, OpenAI's models, particularly \texttt{o3} and \texttt{o4 mini high}, were dominant performers in mathematical reasoning benchmarks. 
Beyond pure reasoning, OpenAI excels in its advanced research and task automation capabilities. The platform offers several powerful tools for information synthesis. A key tool is the \texttt{Deep Research} function, which operates similarly to Google's equivalent but tends to produce more compact and synthesized summaries. Another significant feature, introduced on July 17, 2025, is the \texttt{Agent Mode}. This function transforms the model into a proactive assistant capable of handling complex, multi-step objectives autonomously.

\begin{table}[h!]
	\centering
	\caption{Overview of OpenAI's Subscription Tiers and Models (Until August, 7 2025, prior to the \texttt{GPT 5} consolidation)}
	\label{tab:openai_models}
	\begin{tabular}{p{3.5cm} p{10cm}}
		\toprule
		\multicolumn{2}{c}{\textbf{ChatGPT Pro} --- €229/month} \\
		\cmidrule(lr){1-2}
		\texttt{o3 Pro} & OpenAI's most powerful reasoning model, designed for strong mathematical reasoning. \\
		\multicolumn{2}{l}{\textit{Includes all features from the Plus tier.}}\\
		\midrule
		
		\multicolumn{2}{c}{\textbf{ChatGPT Plus} --- €23/month} \\
		\cmidrule(lr){1-2}
		\texttt{o3} & A strong reasoning model for excelling at complex mathematical reasoning and advanced, in-depth web research. \\
		\texttt{o4 mini high} & A strong reasoning model optimized for strong visual reasoning capabilities. \\
		\texttt{ChatGPT 4.5 preview} & A highly creative model, strongly recommended for writing, brainstorming, and polishing language. \\
		\texttt{Agent Mode} & A proactive assistant that autonomously uses tools to complete complex, multi-step tasks. \\
		\texttt{Deep Research} & After setting up a research plan, explores the web for websites and files for 5-15 minutes. Compact output. 25 queries/month. \\
		\midrule
		
		\multicolumn{2}{c}{\textbf{ChatGPT} --- Free} \\
		\cmidrule(lr){1-2}
		\texttt{4o} & Optimized for fast, conversational responses and high-quality multimodal interaction. Not intended for deep, complex reasoning. \\
		\bottomrule
	\end{tabular}
\end{table}

\subsubsection{The \texttt{GPT 5} Consolidation}
On August 7, 2025, OpenAI fundamentally restructured its product line with the release of \texttt{GPT 5}. This update consolidated the various core reasoning models—namely the \texttt{o}-series and \texttt{ChatGPT 4.5 preview}—into a streamlined, three-tier system. Depending on the task, the \texttt{GPT 5} model automatically selects the most promising agent to fulfill the task. Specialized features such as \texttt{Agent Mode} and \texttt{Deep Research} were not replaced but persist as tools available within the new subscription tiers. The new structure for the core models is as follows:

\begin{itemize}
	\item \textbf{\texttt{GPT 5}:} The new standard model, available in the free tier, replacing the previous \texttt{4o}.
	\item \textbf{\texttt{GPT 5 Thinking}:} The mid-tier subscription model, which consolidates the capabilities of the former \texttt{o3}, \texttt{o4 mini high}, and \texttt{ChatGPT 4.5 preview} models into a single, powerful offering.
	\item \textbf{\texttt{GPT 5 Pro}:} The top-tier premium model that replaced the former \texttt{o3 Pro} as OpenAI's most advanced system for complex reasoning tasks.
\end{itemize}

This consolidation marks a significant shift in OpenAI's strategy, moving from a specialized toolkit of selectable models to a more unified and powerful general-purpose assistant. While this simplifies the user experience, it also means that the specific performance characteristics identified in the \texttt{o}-series benchmarks can no longer be directly attributed to a single, selectable model. The performance of the new \texttt{GPT 5} tiers will require a new wave of independent benchmarking and investigation.

\subsubsection{xAI's Grok Family}
The Grok family of models, developed by xAI, is designed with a distinct philosophy. It is characterized by a streamlined user interface and a powerful, deeply integrated real-time web search capability. While earlier versions of Grok performed solidly in benchmarks, they were not typically considered top-tier contenders. 

This changed with the release of \texttt{Grok 4} on July 10, 2025. As detailed in Subsection \ref{subsec:Competitions in July 2025}, this new iteration quickly ascended to the top of several key mathematical benchmarks, establishing itself as a leading model.

The \texttt{Grok 4} platform is offered in several tiers, as summarized in Table \ref{tab:xai_models}. The standard \texttt{Grok 4} model, a single-agent system, is already highly capable. However, the premium variant, \texttt{Grok-4-Heavy}, employs a more advanced multi-agent architecture. When presented with a complex task, it dynamically spawns multiple independent AI agents that explore different solution pathways in parallel. These agents then collaborate to compare results, verify insights, and synthesize a final, more robust answer. This approach, which significantly increases computational resources at inference time, is designed to enhance performance on particularly challenging problems.

In terms of practical limitations, the platform supports file uploads but offers a smaller context window than its main competitor, with \texttt{Grok 4} handling 128,000 tokens and \texttt{Grok-4-Heavy} managing 256,000. Furthermore, the integrated web search cannot be disabled, which can lead to longer response times during focused document analysis. For tasks requiring a deep, uninterrupted analysis of large, self-contained texts, \texttt{Gemini 2.5 Pro} therefore remains the recommended tool.

However, its development and public integration via the social media platform X (formerly Twitter) have not been without significant controversy. Shortly before the release of \texttt{Grok 4}, a preceding version became notorious for generating a series of antisemitic responses, a widely reported incident that raised serious concerns about the model's alignment and safety protocols \cite{CNNNewsGrok25}. 

For the purpose of this article, this incident serves as a critical, real-world case study. It is perhaps the most vivid public demonstration to date of the inherent risks of large-scale AI models, including unforeseen emergent behavior and the amplification of societal biases present in their training data. It serves as a stark, practical reminder of the principles outlined earlier: AI models are not neutral or infallible agents, and their outputs demand constant and rigorous critical verification from the human user. Therefore, while \texttt{Grok 4} is included in this toolkit for its demonstrated high performance in mathematical reasoning, its use must be accompanied by an acute awareness of these documented vulnerabilities.

\begin{table}[htbp]
	\centering
	\caption{Overview of xAI Models}
	\label{tab:xai_models}
	\begin{tabular}{p{3.5cm} p{3.5cm} p{6cm}}
		\toprule
		\textbf{Model} & \textbf{Cost} & \textbf{Key Strengths} \\
		\midrule
		\texttt{Grok-4-Heavy} & SuperGrok Heavy \$300/month & Maximum reasoning power; high reliability. \\
		\multicolumn{2}{l}{\textit{Includes all features from the SuperGrok tier.}}\\
		\midrule
		\texttt{Grok 4} & SuperGrok \$30/month & Strong reasoning; fast web integration. \\
		\midrule
		\texttt{Grok 3} & Free & basic reasoning. \\
		\bottomrule
	\end{tabular}
\end{table}

\subsubsection{Specialized Auxiliary Tools}
Beyond these general-purpose assistants, mathematicians can leverage specialized tools for specific tasks. A prime example is {DeepL}, a service renowned for its high-quality translation capabilities. Its {DeepL Write} feature, available via the browser or using a desktop application, is particularly effective for polishing academic English, refining sentence structure, and ensuring linguistic precision, making it an invaluable tool in the final stages of manuscript preparation.

\subsubsection{A Note on Coding} \label{subsubsec: A note on coding}
The ability to generate code is a central feature of modern LLMs that enables a powerful workflow for mathematical exploration. The leading platforms discussed—Google's AI Studio, OpenAI's models, and xAI's Grok—all possess robust code generation capabilities, particularly in widely used languages like Python, enabling tasks that range from plotting functions to implementing complex algorithms. While many of these tools offer to execute code directly within their interface, this convenience can come at the cost of transparency. Manually executing the generated code is a direct application of the {Copilot, Not the Pilot} Principle outlined in Subsection \ref{sub: Guiding Principles for AI Collaboration}. This practice ensures full user control, allows for easy modification and debugging, and maintains the reproducibility of any computational work.

Furthermore, for mathematicians with limited programming experience, these AI assistants can lower the barrier to entry for computational work. They can provide step-by-step guidance for setting up a local coding environment, such as installing SageMath within a Windows Subsystem for Linux, effectively removing technical hurdles to experimentation. 

\subsubsection{Data Privacy and Security}
Utilizing cloud-based AI tools requires the transmission of data, including potentially sensitive or unpublished research, to third-party corporate servers. Researchers must be aware that the terms of service for many platforms, particularly free tiers, may permit the use of this data for future model training, effectively incorporating novel ideas into the provider's intellectual property. This practice, combined with the risk of data breaches, necessitates a judgment about what information to share. Thus, adherence to personal and institutional data security policies is essential when handling proprietary or confidential material.

\section{Seven Ways of AI Usage in Mathematical Research }
The preceding analysis of benchmarks, tools, and principles provides the foundation for a practical framework on how to integrate these tools into the mathematical research workflow. This section systematically explores seven fundamental ways in which AI can be integrated into the mathematical research workflow, from initial ideation to final publication. These {Seven Ways} are presented as durable, high-level strategies intended to remain relevant even as the underlying technology evolves.

Within the context of each of these applications, we will pursue a threefold objective. To ensure scientific integrity and to ground our examples in verifiable data, the following illustrations will primarily refer to the models whose performance was documented in our benchmarks, namely \texttt{Gemini 2.5 Pro}, \texttt{Grok 4}, and OpenAI's legacy \texttt{o}-series. The recent release of OpenAI's \texttt{GPT 5} does not invalidate these findings. Instead, it extends these capabilities. The reader is encouraged to map the functionalities discussed for the legacy models to their contemporary counterparts; for instance, the reasoning capabilities attributed to \texttt{o3} and \texttt{o4 mini high} are now included in \texttt{GPT 5 Thinking}. This approach allows us to base our practical advice on established data while guiding the reader on how to apply these principles with the tools available today. Second, we will demonstrate how the abstract Guiding Principles discussed in Subsection \ref{sub: Guiding Principles for AI Collaboration} translate into concrete, necessary practices within these specific workflows. Finally, by examining the strengths and weaknesses of today's tools in these contexts, we aim to develop a deeper understanding of this new technology and how to navigate its potential effectively. As such, this section offers a conceptual map for engaging with AI in mathematical research, both now and in the future.
\subsection{Creativity and Ideation}\label{subsec: Creativity and Ideation}
The research process often begins with a phase of creativity and ideation, a domain where AI can serve as a powerful assistant. While benchmarks like MathArena rightly focus on uncontaminated data to test for genuine reasoning, the vast, \enquote{contaminated} training data of LLMs is precisely what makes them a key resource for generating new ideas. Having been exposed to a volume of mathematical proofs, concepts, and connections far exceeding that of any human, these models can serve as powerful engines for creative exploration.
\subsubsection{Research}
In research, this capability is particularly valuable for formulating new research questions, especially when venturing outside one's primary area of expertise. Because an LLM has been trained on a vast corpus of scientific literature, it possesses a broad awareness of which problems are considered interesting, which are actively being worked on, and what connections exist to adjacent fields—knowledge that a human researcher may not have encountered yet. An effective AI-driven ideation process can leverage this by identifying research questions that are not only tractable but also relevant to an active research community. This includes generating novel examples that satisfy specific criteria or illustrate the application of a theorem. If a proposed research question extends beyond the user's own field, the AI can serve as an initial guide to navigate the new territory, a concept further explored in Subsection \ref{subsec: Interdisciplinarity}. For such creative tasks, models like \texttt{Gemini 2.5 Pro} (with a high {temperature} setting and disabled web search) or \texttt{Grok 4} /\texttt{o3} with web search are particularly well-suited.

The OPC study's finding of an \emph{intuition gap}—a discrepancy between final-answer accuracy and proof validity—suggests that some models possess a form of mathematical intuition. This can be creatively exploited by using a model with high answer accuracy to generate conjectures, which are then rigorously tested. This hybrid approach, balancing machine-generated ideas with human-led verification, is a key aspect of the {Copilot, Not the Pilot} Principle (see Subsection \ref{sub: Guiding Principles for AI Collaboration}).

Furthermore, the best-of-$n$ sampling strategy, discussed in Subsection \ref{subsubsec: The Impact of Best-of-$n$ Sampling on Proof Quality}, is not just a performance booster but also a creative tool. Generating multiple distinct solution paths for a single problem allows a mathematician to explore a wide range of potential strategies. The ability of models to generate code, as noted in Subsection \ref{subsubsec: A note on coding}, further enhances this process by enabling the rapid visualization and empirical testing of abstract conjectures, providing immediate feedback to guide theoretical exploration.

Once a promising idea has been generated through these methods, it can be further investigated using a targeted literature search and literature analysis (see Subsections \ref{subsec: Literature Search}, \ref{subsec: Literatur Analysis}) to ensure a productive workflow. 
\subsubsection{Teaching}
Beyond research, creativity and ideation are also relevant for teaching, a core responsibility for many academic mathematicians. In this domain, the extensive training of LLMs on foundational mathematical textbooks and standard pedagogical material can be a significant asset. The models' familiarity with classic proofs and standard techniques allows them to function as assistants in curriculum development, facilitating the creation of individualized learning materials.

An instructor can utilize these capabilities to construct novel problems that move beyond standard textbook exercises. The process is often iterative: an instructor might begin by using a high {temperature} setting in a model like \texttt{Gemini 2.5 Pro} to generate a wide range of problem scenarios. This approach maximizes output diversity, though it requires the instructor's critical judgment—in line with the {Pilot} principle—to filter out ideas that rely on concepts beyond the course's scope. 

Once a promising direction is identified, or if the instructor has a specific didactic concept in mind, a lower, more deterministic temperature can be used. In this mode, the model can transform a rough idea into a precisely formulated problem, complete with a detailed model solution. The result is a tailored exercise calibrated to the course's syllabus and the students' expected level of knowledge. This strategic use of model parameters, which requires some experimentation, is a direct application of the Experimental Mindset and Model Selection principles outlined in Subsection \ref{sub: Guiding Principles for AI Collaboration}.

\subsection{Literature Search}\label{subsec: Literature Search}
Following the initial ideation phase, a targeted literature search is the logical next step to contextualize and validate new ideas. For a quick overview of a topic, models with strong, integrated web search capabilities, such as \texttt{Grok 4} or OpenAI's \texttt{o3}, are effective. They combine reasoning with real-time information retrieval to provide a preliminary survey of the relevant landscape. It is important to note that querying a model like \texttt{Gemini 2.5 Pro} for literature when its web search is disabled or without providing the source documents can lead to unreliable results, as further discussed in Subsection \ref{subsec: Literatur Analysis}.

For a more in-depth investigation, the specialized \texttt{Deep Research} functions offered by both Google and OpenAI are recommended. These tools perform iterative, automated searches over a period of 5-15 minutes, analyzing websites and documents to synthesize a comprehensive summary. While both function similarly, they exhibit different output styles: Google's tool tends to provide more detailed, expansive summaries, whereas OpenAI's output is typically more compact. A critical advantage of these dedicated search tools is that they provide source citations, enabling the user to directly verify the synthesized information, see the Principle of Critical Verification (Subsection \ref{sub: Guiding Principles for AI Collaboration}). Such specialized systems are the practical result of a dedicated research field known as Mathematical Information Retrieval (MathIR), which addresses the foundational challenges of indexing and searching mathematical content \cite{ZanibbiEtAl25}.

However, users should remain aware that the quality of the output depends on the source material the AI finds. The models may not distinguish between correct and incorrect information within the articles they process and may present flawed results as factual. Once the relevant literature has been identified through this process, the next stage is a detailed literature analysis.

\subsection{Literature Analysis}\label{subsec: Literatur Analysis}
Once a curated list of relevant literature has been compiled, the next step in the research workflow is a detailed analysis of these documents. LLMs offer a suite of capabilities that can significantly accelerate this process, though their use requires an awareness of certain critical limitations.

\subsubsection{Analysis of Provided Documents: Capabilities and Workflow}
Once relevant literature has been identified, LLMs can serve as powerful tools for in-depth analysis. For this task, \texttt{Gemini 2.5 Pro} is highly recommended due to its exceptionally large context window of 1 million tokens. This allows a user to upload multiple research papers or even entire textbooks into a single session, enabling a comprehensive analysis that can be completed in minutes. The model excels at tasks such as summarizing complex arguments, tracking the usage of notation across a lengthy document, or locating specific theorems and cross-references. Unlike traditional text search functions, which are limited to exact string matching, an LLM can understand semantic context and identify concepts even when different variables or phrasings are used. To achieve a good result, it is advisable to set the creativity bar very low and disable the Google search.

A key advantage of this approach is the ability to engage in a mathematical dialogue about the provided material, effectively using the AI as a tireless sparring partner as discussed in Subsection \ref{subsec: Social Aspect}. A researcher can ask for clarifications on difficult passages, discuss the implications of a theorem, or request a comparison of different authors' approaches to a topic.
\subsubsection{A Critical Limitation: Internal Knowledge vs Provided Text}
A strong word of caution is necessary, however, regarding the distinction between analyzing user-provided documents and querying the model's internal knowledge base. When asked to recall the location of a theorem from its training data, a model like \texttt{Gemini 2.5 Pro} (with web search disabled) will often generate a plausible but incorrect citation. It may confidently point to a thematically related source that does not actually contain the specific result in question.

This issue is exacerbated by two fundamental limitations discussed earlier. First, due to its persistent memory (see Subsection \ref{subsubsec: Google DeepMind's Portfolio}), the model becomes strongly anchored to its initial, incorrect belief. Second, this behavior aligns with the reluctance to admit failure mode identified in the OPC study (Subsection \ref{subsubsec: Common Failure Modes in Proof Generation}). The combination of these factors can lead to a critical failure mode: even when the user subsequently provides the correct source document and asks the model to find the information, its prior conviction can be so strong that it fails to locate the passage or, in worse cases, begins to hallucinate, insisting the information is present when it is not. This underscores the critical importance of relying on the model for analysis of explicitly provided texts, rather than for recalling information from its internal, opaque knowledge base.

It is important to note, however, that the analysis capabilities are primarily text-based. While many models can process documents containing images, their ability to accurately interpret complex mathematical graphs, diagrams, or commutative diagrams is often unreliable. The model may generate a plausible-sounding interpretation of a figure regardless of whether it can actually \enquote{read} it, which makes independent verification of any visually-based claims essential.
\subsection{Interdisciplinarity}\label{subsec: Interdisciplinarity}
AI promotes connections across many fields, bridging the gap between mathematics and other sciences, such as physics and engineering, as well as between different subdisciplines of mathematics, such as numerics and analysis. A core function of LLMs in this context is translation. This applies literally, enabling the translation of research papers from languages like Russian or Chinese, and metaphorically, allowing the translation of concepts between the technical \enquote{languages} of different scientific fields.

While this capability makes new areas of study more accessible, it does not replace the need for expert knowledge. On the contrary, when exploring a new domain with an AI, the Principle of Critical Verification (Subsection \ref{sub: Guiding Principles for AI Collaboration}) becomes paramount, and the reliance on the expertise of human collaborators from that field increases. The AI serves as a powerful starting point for exploration, but the validation of complex ideas must ultimately be done in collaboration with human experts.

In such collaborative projects, LLMs can significantly enhance the workflow. A common hurdle between experts from different fields is the difference in technical language and notation. An LLM can act as a \emph{universal translator,} helping each expert to rephrase their concepts for a non-expert audience or to clarify discipline-specific terminology, thereby streamlining the communication process.

Furthermore, an LLM's extensive training in different mathematical fields can be used to find analogies between seemingly unrelated areas. This could lead to the transfer of proof techniques. When an AI suggests such a novel connection, its plausibility can be cross-checked using a multi-model verification strategy. This approach, which involves using a second AI to critique the output of the first, provides a crucial first layer of validation and is discussed in more detail in the context of mathematical reasoning in Subsection \ref{subsec: Mathematical Reasoning}.

Finally, AI provides a powerful bridge between theoretical mathematics and computation. Emulating the success of systems like \texttt{AlphaEvolve} (Section \ref{sec: Landmark Achievements in AI-Driven Mathematical Discovery}), a mathematician can use an AI to translate a theoretical problem into a computational one, generate code to explore it, and use the empirical results to form new conjectures. This process is a direct application of the AI-assisted coding approach discussed in Subsection \ref{subsubsec: A note on coding}. The procedure also applies to interfacing with Computer Algebra Systems (CAS), though the AI's proficiency is limited to systems it has frequently encountered in its training data.
\subsection{Mathematical Reasoning}\label{subsec: Mathematical Reasoning}
The use of AI for direct mathematical reasoning may transforms to its strongest applications. As the benchmarks discussed in Section \ref{sec:From Frontier Models to Publicly Available Tools} demonstrate, leading publicly affordable LLMs like \texttt{Gemini 2.5 Pro}, \texttt{Grok 4}, and OpenAI's \texttt{o3} and \texttt{o4} models possess good reasoning capabilities, not to mention the expensive models such as \texttt{Deep Think}. The dynamic competitive landscape of public LLMs (see Subsection \ref{subsec:Competitions in July 2025}) suggests that the optimal workflow for a complex problem may involve combining the diverse strengths of different models. They can combine known arguments in novel ways to construct proofs for statements they have not encountered before. However, the quality of this reasoning is highly model-dependent, making the selection of a proven, high-performing model a prerequisite for any serious work. This direct application of LLMs to mathematical reasoning is not only a frontier for practitioners but also a key topic within the MathIR community, where researchers investigate how these capabilities can be integrated into systems for search and question answering \cite[Sec.~6.3]{ZanibbiEtAl25}.

\subsubsection{The Interactive Workflow and Principles of Verification}
Engaging with an AI in a reasoning task is an interactive process that requires constant human oversight. Every AI-generated proof or calculation must be critically examined for plausibility and correctness by the researcher, a direct application of the Copilot, Not the Pilot Principle (see Subsection \ref{sub: Guiding Principles for AI Collaboration}). An effective workflow also involves developing an understanding of the specific model in use. Each model exhibits its own systematic biases and failure modes, and learning to anticipate these requires the Experimental Mindset mentioned earlier. A crucial operational detail, particularly for models like \texttt{Gemini 2.5 Pro}, is the management of their persistent memory. Since an incorrect conclusion drawn early in a session may be repeated, it is often essential to start a new chat for a new line of reasoning to ensure a clean logical slate.

Several advanced techniques can further enhance the quality and reliability of AI-assisted reasoning. Given the {self-critique blindness} of LLMs documented in the OPC study (Figure \ref{fig:j2}), a multi-model \enquote{peer review} process might be an effective strategy. A proof sketch generated by one model can be submitted to a different one for verification and refinement, leveraging their strength as external evaluators. Anecdotal experience suggests that while models are adept at finding genuine errors, they can also be over-sceptical; \texttt{Gemini 2.5 Pro}, for instance, sometimes flags correct steps in a human-written proof as potential mistakes. 

Performance can also be significantly boosted by employing best-of-$n$ sampling, as shown by the results for \texttt{o4 mini} (Subsection \ref{subsubsec: The Impact of Best-of-$n$ Sampling on Proof Quality}). Generating multiple solution paths not only increases the probability of finding a correct proof but also provides a range of different strategies to consider. Furthermore, mathematical reasoning can be augmented by code generation (Subsection \ref{subsubsec: A note on coding}); translating a logical step into executable code serves as a powerful computational verification method.

\subsubsection{Reasoning as a Tool for Exploration and Diagnosis}
Beyond proof construction, LLMs serve as powerful interactive tools for exploring and understanding arguments. A researcher can use the AI to quickly explore different lines of reasoning (\enquote{what happens if I apply this theorem here?}) or as a diagnostic tool when stuck on a proof. Presenting a partial argument to the model can often reveal the underlying weakness or missing link in one's own thinking. This ability to engage in a direct, technical dialogue transforms the AI into a constant sparring partner, a concept that leads directly to the social aspects of this new mode of work.
\subsection{Social Aspect}\label{subsec: Social Aspect}
Beyond individual reasoning tasks, the integration of AI into the mathematical workflow introduces significant social dynamics, reshaping how researchers interact with knowledge, their peers, and their students.

\subsubsection{The AI as a Personal Research Sparring Partner}
A key social benefit for researchers is the role of a high-performing AI assistant as a constant sparring partner, available 24/7. It offers an immediate, interactive space to discuss doubts about a proof, explore surprising results from a calculation, or analyze complex literature (as detailed in Subsections \ref{subsec: Mathematical Reasoning} and \ref{subsec: Literatur Analysis}). This constant availability is particularly valuable for reducing the intellectual isolation common in research, especially at smaller institutions. The AI's non-judgmental nature, stemming from its lack of ego, creates a safe environment to test nascent or unconventional ideas without the hesitation one might feel with a human colleague.

\subsubsection{Enhancing Human and Interdisciplinary Collaboration}
AI can also significantly enhance collaboration between researchers, particularly in interdisciplinary projects (see Subsection \ref{subsec: Interdisciplinarity}). The accessibility of powerful free tools like Google's AI Studio creates a shared technological baseline, allowing researchers from different institutions to work with comparable computational support. In direct collaboration, an LLM can act as a neutral \enquote{third-party arbiter.} When experts disagree on a technical point, they can submit it to a trusted model for an unbiased assessment, leveraging the near-human accuracy of top models as evaluators (as shown in Figure \ref{fig:j1}) to resolve disputes efficiently.

\subsubsection{Implications for Teaching and Learning}
The social impact on education is twofold. For instructors, the ability to automate the creation of routine exercises, as discussed in Subsection \ref{subsec: Creativity and Ideation}, can free up time for other pedagogical activities such as mentoring or leading seminars. For students, the AI can serve as a personal tutor and sparring partner, analogous to its role in research. It provides a space for one-on-one support where they can explore concepts and ask questions at their own pace.

However, this application requires particular care. Unlike experienced researchers, students are often less equipped to distinguish between correct reasoning and plausible-sounding falsehoods. This elevates the importance of using a highly reliable and factually grounded reasoning model for educational purposes, and it underscores the need for students to understand the non-human nature of the AI, as outlined in our third guiding principle. This functionality is being developed further with specialized features, such as the \texttt{Study Mode} introduced by OpenAi in July 2025, which offers interactive quizzes, hints and adapts its responses to the user's knowledge level. While the intended outcome of such tools is to foster independent learning, their effective and safe use hinges on a critical and mindful approach from both the student and the instructor.
\subsubsection{Warning: The Non-Human Nature of the AI Partner}
It is crucial to remember that the AI is not human. Attributing human-like qualities to it is a fundamental error that can lead to flawed interactions. Two key differences are particularly important. First, LLMs are designed to generate convincing and fluent statements, which can be deceptive; a well-phrased argument is not necessarily a correct one, demanding a constant critical stance towards all outputs. Second, as demonstrated by the OPC study, LLMs exhibit a reluctance to admit failure (see Subsection \ref{subsubsec: Common Failure Modes in Proof Generation}). This means that, unlike a human, an AI does not easily discard a disproven belief. A flawed conclusion from earlier in a conversation can resurface, a behavior particularly noticeable in models with large, persistent memory contexts.
\subsection{Writing}
The writing phase represents the culmination of the research process, where ideas are formalized and communicated. AI tools can assist in nearly every aspect of this stage, from structuring initial drafts to preparing the final manuscript for dissemination.

\subsubsection{Structuring and Refining the Core Argument}
A primary challenge in academic writing is transforming non-linear research thoughts into a linear, logical argument. An AI can assist in this crucial translation step by structuring rough notes, partial proofs, and intuitive ideas from the \emph{Creativity} and \emph{Reasoning stages} into a coherent draft, complete with formal definitions, theorems, and proofs. Once a draft is established, maintaining consistency in notation, terminology, and argumentation across a lengthy paper is critical. Here, LLMs with a large context window, like \texttt{Gemini 2.5 Pro} on low temperature, excel. A researcher can upload their entire manuscript and task the AI with acting as a \emph{consistency checker,} identifying inconsistencies that a human author might overlook.

For refining the language, specialized tools like DeepL Write are highly effective for polishing sentences and ensuring linguistic precision. General-purpose models such as \texttt{ChatGPT 4.5 preview} or \texttt{Gemini 2.5 Pro} (with a low temperature setting) are also strong at improving prose, but require precise prompting to avoid unintentionally altering the mathematical content. The quality of the stylistic output is highly dependent on the user's prompting skill, as emphasized in the guiding principles.

\subsubsection{Auxiliary Tasks and Ethical Considerations}
Beyond drafting the core text, AI can streamline several auxiliary writing tasks. This includes finding relevant references for a claim or helping to maintain a clear and consistently formatted bibliography. Furthermore, LLMs are exceptionally adept at distilling a complex paper into a concise and compelling abstract. A researcher can provide the full text and ask for several abstract versions, which can then be refined. After a paper is written, the work often needs to be presented to different audiences; an LLM can efficiently adapt the dense, technical language of the paper into various styles and levels of detail, aiding in the effective distribution of research.

This phase also makes the ethical considerations of AI usage (discussed in Subsection \ref{subsec: Navigating the Ethical Landscape: Authorship, Plagiarism, and Scientific Responsibility}) most concrete. While using AI for tasks like language polishing is widely accepted, generating core parts of the text requires careful handling. This reinforces the principle of \enquote{intellectual ownership}: by the time of submission, the author must have fully verified and internalized the text, ensuring it is authentically their own. A brief mention of the tools used in the acknowledgements is becoming a best practice for transparency.
\subsection{Ethical Considerations: Authorship and Responsibility}\label{subsec: Navigating the Ethical Landscape: Authorship, Plagiarism, and Scientific Responsibility}

The integration of AI into the core of the research process, particularly in mathematical reasoning, raises significant and methodological questions. If an AI generates a key step in a proof, who is the author? Does using such a tool constitute plagiarism? And what is our responsibility as researchers in a world where other scientific fields, as seen with \texttt{AlphaFold}'s impact on biology, are achieving Nobel Prize-worthy breakthroughs by embracing these new technologies?

These are not questions with simple answers, but navigating them requires a clear framework. The concern of plagiarism stems from viewing the AI as an independent author. However, a more accurate model is that of a highly sophisticated, interactive instrument. Just as a mathematician using a computer algebra system is not plagiarizing the software's creators, a researcher using an LLM is directing a tool. The intellectual ownership and responsibility for the final output—its verification, refinement, and contextualization—remain entirely with the human pilot, as emphasized in Subsection \ref{sub: Guiding Principles for AI Collaboration}. This necessity for human oversight is not merely a matter of current technological immaturity. Indeed, theoretical arguments suggest that this human role may be a permanent feature of mathematical research. As explored by Dean and Naibo, many significant mathematical problems possess a logical structure that makes them inherently resistant to purely algorithmic resolution \cite{SiliconReckoner DeanNaibo25}. This perspective implies that the conceptual leaps and rigorous validation required to solve such problems will continue to demand human intellect, reinforcing the pilot's ultimate responsibility for the final result.

The question of whether \textit{not} to use these tools is perhaps even more pressing. Mathematical progress has often been accelerated by the integration of new computational instruments, from the first mechanical calculators to modern computer algebra systems. To ignore tools that are proven to enhance our reasoning and creativity would mean neglecting potential paths to discovery. The real question is not whether today's models are essential for every problem, but how to responsibly engage with a technology that is improving so quickly. Therefore, as long as the process remains transparent and the final work is rigorously verified by the human expert, using AI to augment research is a natural continuation of the scientific tradition of tool adoption. The growing accessibility of these capable models strengthens their potential to act as a common instrument for discovery. As these tools become more integrated into the research workflow, establishing clear community standards for their transparent acknowledgment will be essential.

\section{Conclusion and Future Outlook}

The integration of artificial intelligence into mathematics marks a crucial moment, shifting the paradigm of how research is conducted. This article has explored this evolving landscape, revealing a complex duality of augmentation and automation. While systems like \texttt{AlphaEvolve} point towards a future with increasing automation of specific research tasks, the current and foreseeable impact of AI for the working mathematician is primarily one of augmentation. The model of the AI as a {copilot} acknowledges this tension: it enhances the researcher's capabilities, but places the ultimate responsibility for strategic direction, critical judgment, and intellectual ownership firmly with the human {pilot.}

A clear picture of the capabilities of the current generation of publicly affordable, state-of-the-art models emerges when viewed across the full spectrum of mathematical competitions. They achieve high accuracy on pre-Olympiad problems, demonstrating a strong command of multi-step logical reasoning. Their performance on the exceptionally difficult Olympiads—tasks that would challenge many professional mathematicians—provides a realistic measure of their current frontier. However, the stark contrast between this measured frontier for affordable models and the gold-medal success of the premium-tier model \texttt{Gemini Deep Think} suggests that this is not a permanent technological barrier, but rather a sign of a rapidly maturing field. The mathematical community must therefore prepare for a future of continuously improving AI reasoning capabilities.

To navigate this dynamic environment, this article has proposed a durable framework for practice. The {Seven Ways} of AI Usage offer a stable set of strategies for integrating these tools into the research workflow, designed to remain relevant even as the specific models evolve. Central to this framework is a renewed emphasis on the researcher's role as the ultimate judge of correctness. As AI tools can produce convincing but flawed arguments, the core scientific standard of rigorous, independent verification becomes quite important. Engaging with AI under this principle of critical oversight transforms the research process: it empowers the mathematician to explore more ambitious questions, to test connections between disparate mathematical fields, and to verify lines of reasoning that might otherwise remain unexplored.

Looking forward, the future of AI in mathematics will likely involve a deeper integration of natural language models with formal proof systems and the development of more specialized AI agents. As these tools evolve, so too will the research practices that effectively leverage them, continually reshaping the process of mathematical discovery. This implies a necessary evolution in mathematical training. Future mathematicians will need to be skilled not only in their specific domain but also in the art of human-AI interaction—including strategic prompting, critical output evaluation, and the ethical navigation of these tools. Integrating AI skills into graduate education will be a key challenge and opportunity for the years to come.
\section*{Acknowledgements}
I would like to express my sincere gratitude to my brother, Jörn Henkel, for introducing me to the latest AI tools, and to my supervisor, Ilka Agricola, for encouraging me to write this article. 

The creation of this text itself served as a practical application of the principles discussed herein; various AI instruments were utilized throughout the research and writing process, most notably Google's \texttt{Gemini 2.5 Pro} for ideation and analysis, \texttt{o3} and \texttt{Gemini Deep Research} for literature search and DeepL as well as \texttt{ChatGPT 4.5 preview} for language refinement.

\section*{Conflict of Interest}
The author declares that he has no conflict of interest.		
	
%----------------------------------------------------	
%----------------------------------------------------	
%----------------------------------------------------

\end{document}